\documentclass[a4paper]{article}

\usepackage{graphicx}
\usepackage{float}
\usepackage{subfigure}
\usepackage{amsfonts}
\usepackage{amsmath}
\usepackage{amsthm}
\usepackage{amssymb}
\usepackage{bm}
\usepackage{authblk}
\usepackage{color}
\usepackage[T1]{fontenc}
\usepackage[utf8]{inputenc}
\usepackage{comment}
\usepackage{mathtools}
\usepackage{hyperref}
\usepackage{xcolor}
\usepackage{url}

\usepackage[backend=bibtex]{biblatex}
\usepackage{lineno}

\newenvironment{system}%
	{\left\lbrace \begin{array}{@{} l @{} }}%
	{ \end{array}\right.}

\newtheorem{pro}{Proposition}

\newtheorem{thm}{Theorem}
\newtheorem{defi}{Definition}

\newcommand{\bbR}{\mathbb{R}}

\newcommand{\bbP}{\mathbb{P}}

\newcommand{\calD}{\mathcal{D}}
\newcommand{\calE}{\mathcal{E}}
\newcommand{\calI}{\mathcal{I}}
\newcommand{\calR}{\mathcal{R}}
\newcommand{\bbC}{\mathbb{C}}
\newcommand{\calO}{\mathcal{O}}

\newcommand{\calC}{\mathcal{C}}
\newcommand{\calX}{\mathcal{X}}
\newcommand{\bbSigma}{\bm{\Sigma}}
\newcommand{\bOmega}{\bar{\Omega}}

\DeclareMathOperator{\ctg}{ctg}
\DeclareMathOperator{\cst}{cst}
\DeclareMathOperator{\im}{Im}
\DeclareMathOperator{\erf}{erf}
\DeclareMathOperator{\p}{\phi}

\newcommand{\norm}[1]{\left\lVert#1\right\rVert}

\title{A linear barycentric rational interpolant on starlike domains}

\author[a]{Jean-Paul Berrut}
\affil[a]{Department of Mathematics, Université de Fribourg, Fribourg, Switzerland, \\ \texttt{jean-paul.berrut@unifr.ch}}

\author[b]{Giacomo Elefante}
\affil[a]{University of Padova, Department of Mathematics \lq\lq Tullio Levi-Civita\rq\rq, Italy, \\ \texttt{giacomo.elefante@unipd.it}}

\date{}

\begin{document}

\maketitle

\begin{abstract}
When an approximant is accurate on the interval, it is only natural to try to extend
it to several-dimensional domains. In the present article, we make use of the fact that
linear rational barycentric interpolants converge rapidly toward analytic and
several times differentiable
functions to interpolate on two-dimensional starlike domains parametrized in polar coordinates.
In radial direction, we engage interpolants at \textit{conformally shifted Chebyshev
nodes}, which converge exponentially toward analytic functions.
In circular direction, we deploy linear rational trigonometric barycentric interpolants,
which converge similarly rapidly for periodic functions,
but now for \textit{conformally shifted equispaced nodes}.
We introduce a variant of a tensor-product interpolant of the above two schemes
and prove that it converges exponentially
for two-dimensional analytic functions---up to a logarithmic factor---and
with an order limited only by the order of differentiability for real functions (provided that the boundary enjoys the same order of differentiability). Numerical examples confirm that the shifts permit to reach a much higher accuracy with significantly fewer nodes, a property which is especially important in several dimensions.

\end{abstract}



\vspace{1cm}

Interpolating a function in several dimensions is a fundamental research topic in applied mathematics, due to its many applications in engineering and many other fields; one must only recall finite element and pseudospectrals methods for partial differential equations.
In the present work, we shall merely be interested in infinitely smooth interpolants. For the abundant literature on splines, the reader may consult \cite{Splines}, among the many books on the subject.

Classical multivariate rational interpolation has been widely studied  \cite{CuytVerdonkMultivariate,CuytErrorFormulae}. As in the unidimensional case, it has the ability to better deal with singularities of the function than polynomials \cite{DeVoreMultivariate}. Unfortunately, rational approximation can be numerically fragile to compute and is inclined to spurious singularities. Many different approaches have been introduced to avoid these drawbacks, such as algorithms based on the singular value decomposition \cite{Aerospace} or on a Stieltjes procedure and an optimization formulation \cite{AustinMultivariate}.

In this work we aim to construct smooth surfaces by interpolating over non-trivial two-dimensional domains, e.g., domains that are neither rectangles nor disks.
Since the number of interpolation points grows very rapidly in several dimensions
(\lq\lq curse of dimensionality\rq\rq), we want to use fast converging approximants such as the polynomial interpolant at Chebyshev points. However, even the latter is inaccurate for functions with large gradients, as documented on the first line of each of the tables in \cite{BerMitPROC}.
Put in another way, a large number of nodes in every direction (quickly more than 1'000) is necessary for an accurate approximation, which is difficult in several dimensions on a classical personal computer.

It has been known for a long time that in such situations one should concentrate nodes in the vicinity of large gradients (in the context of the solution of time evolution PDEs one speaks of moving meshes).
To maintain the rapid convergence of the infinitely smooth interpolant, one uses for that purpose
infinitely differentiable or even analytic point shifts \cite{BaylissTurkel}. Such shifts may be performed in (at least) two ways \cite{BalBerMit}.
One is to change the variable before using the interpolating polynomial; however, the inner derivatives involved in the chain rule rapidly complicate the formulas when the order of derivation increases.
Here we shall rather use linear rational barycentric interpolation, whose derivatives may be derived from simple formulas discovered by Schneider and Werner \cite{SchneiderWerner} and retain the barycentric nature of the interpolant.
The effect of the point shifts on the accuracy is spectacular (compare the fifth line with the first in the tables of \cite{BerMitPROC}).

For simplicity and cheapness, the method we suggest works with a representation of the domain in
polar coordinates and is therefore limited to starlike domains.
It transplants the problem into a disk and basically consists in a tensor product of a linear rational barycentric interpolant at conformally shifted Chebyshev nodes on the segments in the radial direction with a linear rational trigonometric barycentric interpolant at conformally shifted equispaced nodes in the circular direction.

To retain the rapid convergence of the one-dimensional interpolants, one should obtain these when inserting a constant value of the other variable: for instance, when the radial variable is constant, one should obtain a trigonometric barycentric one-dimensional interpolant in the circular variable.
To guarantee this in a simple and cheap way, we represent the domain through a homothety of its boundary, i.e., by multiplying the polar representation of the latter by all numbers between 0 and 1, and the interpolation points by the intersections of such curves with radii from the origin.
The concentration of nodes at locations of large gradients is obtained by concentrating homothetic curves as well as radii in the vicinity of such locations.

In the next section, we review the one-dimensional linear rational barycentric interpolants which make up the basis of the suggested method, together with the corresponding convergence results (high degree algebraic convergence for several times differentiable functions, exponential convergence for analytic ones). Section 2 recalls the tensor product interpolation of (linear) rational barycentric interpolants on disks, describes the special case of the interpolants used in our method and studies the Lebesgue constant of the resulting two-dimensional interpolant, which is interesting for its own sake and as an important ingredient of the convergence theorems of the final scheme.
Section 3 introduces the formula for the interpolant at shifted nodes on starlike domains and studies its Lebesgue constant, which arises in the error estimates of Section 4. Section 5 demonstrates the efficiency of the method with numerical examples. The paper concludes with
some general remarks on the scheme. 

\section{Preliminaries}

Let us consider $n+1$ distinct nodes in $[a,b]$, which we denote by $x_i$, $i=0,\dots,n$, with their corresponding data $f_i$. In case we interpolate a function $f$ we define the $f_i$ as $f_i:=f(x_i)$.

Then, it is well known (see, e.g., \cite[p.\ 238]{EquiWeights}) that the unique polynomial of degree at most $n$ interpolating the $f_i$ may be written in the barycentric form

\begin{equation}
    p_n(x) = \frac{ \sum_{i=0}^n \frac{\lambda_i}{x-x_i} f_i  }{\sum_{i=0}^n \frac{\lambda_i}{x-x_i}},
\end{equation} 
with the so-called weights $\lambda_i$ defined by
$$ \lambda_i = \prod_{\substack{j=0 \\ j\neq i}}^n \frac{1}{x_i-x_j}. $$

In case the nodes are the Chebyshev points of the second kind, i.e.,
$$ x_i =- \cos\left( \frac{i \pi}{n} \right), \quad i=0,\dots,n,$$
the weights are (up to a constant) \cite[p.\ 252]{EquiWeights}
\begin{equation} \label{ChebWeights}
    \lambda_i = (-1)^i \delta_i, 
\end{equation}
with $\delta_i$ equal to $1/2$ for the first and the last nodes and to $1$ for the others.
The interpolant therefore becomes
\begin{equation} \label{ChebyInt}
    p_n(x) = \dfrac{ \sum_{i=0}^n {\vphantom{\sum}}'' \frac{(-1)^i}{x-x_i} f_i }{\sum_{i=0}^n {\vphantom{\sum}}'' \frac{(-1)^i}{x-x_i}},
\end{equation}
where the double prime means that the first and the last terms of the sum are halved.

With these nodes, $p_n$ enjoys excellent convergence properties when the interpolated function is smooth (see, e.g., \cite[p.\ 53]{Trefethen:ATAP}). To lighten the notation we consider the infinity norm whenever we do not explicit it.

\begin{thm} \label{ChebyErrorJump}
Let $\nu\geq 1$ be an integer and let $f$ and its derivatives up to $f^{(\nu-1)}$ be absolutely continuous on $[-1,1]$ and the $\nu^{\text{th}}$-derivative be of bounded variation, with $V$ its total variation. Then for any $n\geq \nu$ the Chebyshev interpolant satisfies
    $$ \norm{f-p_n} \leq \frac{4V}{\pi \nu(n-\nu)^\nu}.$$ 
\end{thm}

In the case of functions analytic not only in the interval $[-1,1]$ but also in an ellipse with foci $-1$ and $1$, called the \textit{Bernstein ellipse} and denoted by $E_\sigma$, where $\sigma$ is the sum of the lengths of the semiminor and the semimajor axes, the convergence is geometric (see e.g. \cite[p.\ 143]{RivlinCheb}).

\begin{thm} \label{ChebError}
Let $f$ be a function analytic inside and on an ellipse $E_\sigma$. Then for each $n\geq 0$ the polynomial interpolant between Chebyshev points of the first or the second kind satisfies 
    $$ \norm{f-p_n}  \leq \frac{4M}{\sigma^{n}(\sigma-1)} ,$$
where $M=\max_{z\in E_\sigma} |f(z)|$.
\end{thm}

For other nodes than Chebyshev points, the expression \eqref{ChebyInt} becomes a linear rational interpolant \cite{BerInt1}.
Baltensperger et al. proved in \cite{BerBalNoel} that the exponential convergence then still holds for Chebyshev nodes shifted with a conformal map, say $g$. In fact, in this case we have the following theorem for the interpolant
\begin{equation} \label{BarIntTransf}
    r_n(x) = \frac{ \sum_{i=0}^n {\vphantom{\sum}}'' \frac{(-1)^i}{g(x)-g(x_i)} f(g(x_i)) }{\sum_{i=0}^n {\vphantom{\sum}}'' \frac{(-1)^i}{g(x)-g(x_i)}}.
\end{equation}

\begin{thm} \label{TransfConv}
Let $\calD_1$, $\calD_2$ be two domains of $\bbC$ containing $J=[-1,1]$, respectively $I(\subset\bbR)$, let $g$ be a conformal map $\calD_1\to\calD_2$ such that $g(J)=I$, and $f$ be a function $\calD_2\to\bbC$ such that the composition $f\circ g:\calD_1\to\bbC$ is analytic inside and on an ellipse $E_\sigma$ $(\subset \calD_1)$, $\sigma>1$, with foci $-1$ and $1$ and with the sum of the semiminor and the semimajor axes equal to $\sigma$. Let $r_n$ be the interpolant \eqref{BarIntTransf} interpolating $f$ at the transformed Chebyshev nodes $y_i=g(x_i)$, $i=0,\dots,n$. Then
$$ \norm{f-r_n}=\calO(\sigma^{-n}).$$
\end{thm}

Note that, if we write $y=g(x)$ and $y_i=g(x_i)$ for $i=0,\dots,n$, then $r_n$ can be seen as

\begin{equation*} 
    \tilde{r}_n(y) = \frac{ \sum_{i=0}^n {\vphantom{\sum}}'' \frac{(-1)^i}{y-y_i} f(y_i) }{\sum_{i=0}^n {\vphantom{\sum}}'' \frac{(-1)^i}{y-y_i}},
\end{equation*}
which corresponds to Berrut's second interpolant \cite{BaryHormann} $R_1$, introduced in \cite{BerInt1}, at the nodes $y_i$.

\vspace{6pt}
We now turn to periodic interpolants. For the equidistant nodes 
\begin{equation}\label{eq:equinodes}
    \theta_i=i\cdot \frac{2\pi}{n},\qquad \text{ for } i=0,\dots,n-1,
\end{equation} 
in $[0,2\pi)$, there exists a unique balanced trigonometric polynomial that interpolates the data $f_i$, $i=0,\dots,n-1$, at the nodes, and it can be written in the simple barycentric form

\begin{equation} \label{BaryTrigEqui}
    T_n(\theta) = \frac{ \sum_{i=0}^{n-1} (-1)^i \cst\left( \frac{\theta-\theta_i}{2}  \right) f_i}{\sum_{i=0}^{n-1} (-1)^i \cst\left( \frac{\theta-\theta_i}{2}  \right)},
\end{equation} 
where the function $\cst$ is 
\begin{equation} \label{cst}
    \cst(\theta) \coloneqq \begin{system} \csc(\theta), \quad \text{ if } n \text{ is odd,} \\
\ctg(\theta), \quad \text{ if } n \text{ is even}.
\end{system}
\end{equation}

Let us define $\Pi$ as the class of periodic complex valued functions on $\bbR$ with period $2\pi$. Then, it is possible to prove convergence theorems for $T_n$ 
which depend on the smoothness of the function $f$, as with Chebyshev interpolation (see \cite[p. 365]{EquiWeights}, and \cite{PeriodicChebfun} 
for the first theorem and \cite{ConvFourier} for the second).

\begin{thm} \label{TrigoErrorJump}
Let $f\in \Pi$ be such that it has at most simple jump discontinuities in the $\mu^{\textrm{th}}$-derivative which is of bounded variation $V$; then
$$ | f(\theta)-T_n(\theta)|\leq \frac{V}{\pi m^\mu}\left(\frac{1}{m}+\frac{2}{\mu} \right),$$
with $m:=\lfloor n/2 \rfloor$.
\end{thm}

\begin{thm} \label{TrigoError}
Let $f\in\Pi$ be analytic in the strip $S_a:=\{\eta : |\im(\eta)|\leq a \}$, where $a>0$, and $|f(\theta)|\leq M$. Then
$$ | f(\theta)-T_n(\theta)|\leq  2M \cot \left( \frac{a}{2}\right) e^{-an}.$$
\end{thm}

The first author proposed in \cite{BerInt1} to use $T_n$ for nodes other than equidistant points; then \eqref{BaryTrigEqui} becomes a linear rational trigonometric interpolant, which we denote by $t_n \equiv t_n[f]$.

In the wake of the barycentric rational interpolant at the Chebyshev nodes, Baltensperger \cite{Baltensperger} showed in 2002 that, if we consider as nodes some conformally shifted equidistant points, then the interpolant between these shifted nodes retains the convergence behaviour of the one at equidistant nodes. In fact the following holds.

\begin{thm} \label{BalThmConf}
Let $g\in\Pi$ be a conformal map such that $g(I)=I$ with $I=[0,2\pi]$ and such that the function 
$$ w(\p,\theta):= \frac{ \cst\Big( \frac{g(\p)-g(\theta)}{2} \Big) }{ \cst\Big( \frac{\p-\theta}{2} \Big) } $$
is bounded and analytic in $S_{a_1}\times S_{a_1}$ for an $a_1>0$. Let $f$ be a function such that $f\circ g \in \Pi$. Let 
\begin{equation}
  t_n(\p)\equiv t_n[f\circ g] (\theta)= \frac{ \sum_{i=0}^{n-1} (-1)^i \cst\Big( \frac{g(\theta)-g(\theta_i)}{2}  \Big) f(g(\theta_i))}{\sum_{i=0}^{n-1} (-1)^i \cst\Big( \frac{g(\theta)-g(\theta_i)}{2}  \Big)}
\end{equation}
be the rational function generalizing \eqref{BaryTrigEqui} between the nodes $\p_i:=g(\theta_i)$, $i=0,\dots,n-1$, where the $\theta_i$'s are the equidistant nodes \eqref{eq:equinodes}. Then for every $\p\in I$ we have the following.
\begin{itemize}
    \item If $f\circ g$ has simple jump discontinuities in the $\mu^{\textrm{th}}$-derivative, then 
            $$ |f(\p)-t_n(\p)|= \calO(n^{-\mu}).$$
    \item If $f \circ g$ is analytic in a strip $S_{a_2}$ with $a_2\geq a_1>0$ and $f$ is bounded, then
            $$ |f(\p)-t_n(\p)|= \calO(\sigma^{-n}) \quad \text{ with } \sigma:=e^{a_2}.$$
\end{itemize}
\end{thm}

\section{The tensor product interpolant and its Lebesgue constant} \label{Sec2}

Let the fixed integer $d\geq 2$ be the dimension of the domain of the functions to be interpolated. Then, for each, $\ell=1,\dots d,$ consider a set of interpolating functions $\{b_i^{(\ell)}(x)\}$ over $n_\ell+1$ nodes $x_i$ in an interval $I_\ell$ such that
$$ b_i^{(\ell)}(x_k) = \delta_{i,k} $$
and forming the one-dimensional linear interpolant 
$$ \calI_{n_\ell}^{\ell} f(x) = \sum_{i=0}^{n_\ell} f(x_i) b_i^{(\ell)}(x). $$
To construct a $d$-dimensional tensor product interpolation operator of a multivariate function 
$$ f:=f(\xi_1,\dots,\xi_d), \quad f:I^d\to \bbR$$ 
in the box $I^d = I_1 \times \dots \times I_d$, we consider the operator
$$\calI: \calC(I^d)\to \calC(I^d),\qquad \calI = \calI_{n_1}^1 \otimes \dots \otimes \calI_{n_d}^d,$$
where $\calI_{n_j}^j$ is a one dimensional linear interpolation operator over a set of nodes $x_i^{(\ell)}\in I_\ell$, with $i=0,\dots n_\ell$. Then, the tensor product interpolant is of the form (see e.g. \cite[p.\ 30]{TensorNumericalMethods}) 
\begin{equation}
    (\calI f)(\xi_1,\dots,\xi_d) = \sum_{i_1=0}^{n_1}\cdots \sum_{i_d=0}^{n_d} f(x_{i_1}^{(1)},\dots,x_{i_d}^{(d)}) b_{i_1}^{(1)}(\xi_1)\cdots b_{i_d}^{(d)}(\xi_d).
\end{equation}

We will now focus on the case $d=2$ and the domain 
\begin{equation} \label{box}
    B=[0,2]\times [0,2\pi).
\end{equation} 
Starting with the first interval, we let $\calX_n=\{x_0,\dots,x_n\}$ be the set of distinct nodes, in $[0,2]$, and use the basis functions
$$ b_i^{(1)}(x) = \frac{ \frac{(-1)^i \delta_i \eta_i }{x-x_i} }{ \sum_{j=0}^{n_1} \frac{(-1)^j \delta_j \eta_j}{x-x_j}}, $$
in the radial direction, where $\delta_i$ is as \eqref{ChebWeights}  and $\eta_i$ defined by
\begin{equation} \label{weightsR1}
    \eta_i := \begin{system} \sqrt{1-(x_i - 1)^2}, \quad\;\, 0,2\notin \calX_n,\\
\sqrt{\frac{1+(x_i - 1)^2}{2}}, \qquad\quad 0\notin \calX_n, 2\in\calX_n, \\
\sqrt{\frac{1-(x_i - 1)^2}{2}},\qquad\quad 0\in\calX_n, 2\notin\calX_n, \\
1, \qquad\qquad\qquad\quad\,\,\, 0,2\in\calX_n.
\end{system}
\end{equation}
($b^{(1)}_i(x)$ is the transplantation of $b^{(2)}_i(\phi)$, defined below, to the interval $[0, 1]$ through the change of variable $\phi = \arccos x$; the formula for $\eta_i$ is on $[0, 2]$ that of the weights of the polynomial interpolating at the four kinds of Chebyshev points in this interval, see \cite{BerInt1}).
In the angular direction, we take $n_2$ distinct nodes $\phi_0,\dots,\phi_{n_2-1}\in[0,2\pi)$ and we consider the periodic basis
$$ b_i^{(2)}(\p) = \frac{(-1)^i \cst \left( \frac{\p-\p_i}{2} \right) }{\sum_{j=0}^{n_2-1} (-1)^j \cst \left( \frac{\p-\p_j}{2} \right)}.$$
On two-dimensional domains this leads to the tensor product 
$$ \calI = R_1 \otimes t_n $$
of Berrut's second interpolant $R_1$ and the trigonometric interpolant $t_n$ introduced in \cite{BerInt1}, and is an operator from the space of continuous functions in $B$ to $\calR_n \otimes \bbSigma_n$, where $\calR_n$ is the space of rational interpolants with the fixed denominator $\sum_{j=0}^{n_1} ( (-1)^j \delta_j \eta_j ) /(x-x_j)$ and $\bbSigma_n$ is the space of rational trigonometric functions of degree $\lfloor \frac{n_2}{2}\rfloor$ with the same fixed denominator $\sum_{j=0}^{n_2-1} (-1)^j \cst\left( \frac{\p-\p_j}{2}\right)$.

Note that, since the $b_i^{(2)}$ are periodic, this interpolant can be seen as one in polar coordinates in the disk
\begin{equation} \label{DiskRad2}
    E = \{ x\in \bbR^2 : \norm{x}_2\leq 2\}. 
\end{equation} 

\begin{figure}
    \centering
    \graphicspath{ {Images/}}
    \subfigure[]
        {\includegraphics[height=1.74in]{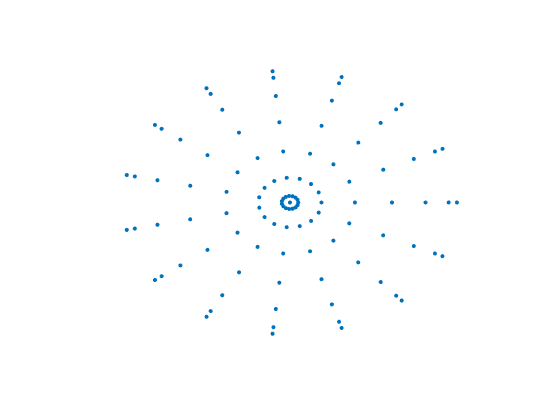}}
     ~ 
    \subfigure[]
        {\includegraphics[height=1.74in]{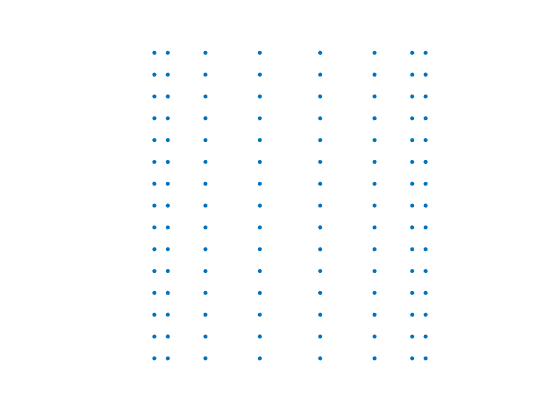}}
    \caption{An example of grid with $n_1=7$ and $n_2=15$, in the disk and in the rectangle $B$ }
    \label{ChebEquiGRID}
\end{figure}

To take advantage of the rapid unidimensional convergence and good conditioning of those interpolants, we will consider as nodes the direct product of Chebyshev points of the second kind in $[0,2]$,
$$ r_i =  -\cos\left( \frac{i\, \pi}{n_1} \right) +1  ,\qquad i=0,\dots,n_1,$$
and equidistant nodes in $[0,2\pi)$,
$$ \theta_j = j\cdot \frac{2\,\pi}{n_2}, \qquad j=0\dots,n_2-1 .$$
An example of such a direct product is shown in Figure \ref{ChebEquiGRID}. Then, the interpolant will be of the form 
\begin{equation} \label{2dBaryInterp} 
    \calI [f] (r,\theta) = \frac{ \sum_{i=0}^{n_1} {\vphantom{\sum}}'' \sum_{j=0}^{n_2-1} \frac{(-1)^{i+j}}{r-r_i} \cst\left(\frac{\theta-\theta_j}{2} \right) f(r_i,\theta_j)    }{\sum_{i=0}^{n_1}{\vphantom{\sum}}'' \sum_{j=0}^{n_2-1} \frac{(-1)^{i+j}}{r-r_i} \cst\left(\frac{\theta-\theta_j}{2} \right) }.
\end{equation}

Note that using these bases allows us to consider conformally shifted nodes in the two variables and thus the interpolant 
\begin{equation} \label{2dConfBaryInterp}
    \calI [f] (r,\theta) = \frac{ \sum_{i=0}^{n_1} {\vphantom{\sum}}'' \sum_{j=0}^{n_2-1} \frac{(-1)^{i+j}}{g_1(r)-g_1(r_i)} \cst\left(\frac{g_2(\theta)-g_2(\theta_j)}{2} \right) f(g_1(r_i),g_2(\theta_j))    }{\sum_{i=0}^{n_1} {\vphantom{\sum}}'' \sum_{j=0}^{n_2-1} \frac{(-1)^{i+j}}{g_1(r)-g_1(r_i)} \cst\left(\frac{g_2(\theta)-g_2(\theta_j)}{2} \right) }
\end{equation}
with the two unidimensional conformal maps 
$$ g_1:[0,2]\to[0,2]\quad \text{ and} \quad g_2:[0,2\pi)\to[0,2\pi).$$

\begin{figure}[H]
    \centering
    \graphicspath{ {Images/}}
    \subfigure[]
        {\includegraphics[height=1.74in]{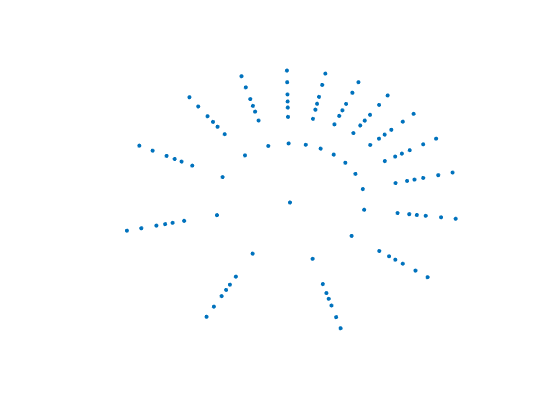}}
     ~ 
    \subfigure[]
        {\includegraphics[height=1.74in]{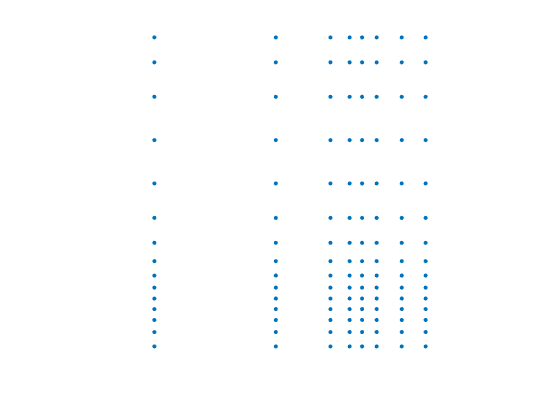}}
    \caption{An example of grid with shifted nodes clustered around $(0.75,\pi/3)$ in the circle and the rectangle $B$ by using the Bayliss-Turkel map for $r$ and the map introduced in \cite{BerElef} for $\theta$ and with $n_1=7$ and $n_2=15$  }
    \label{ChebEquiGRID_Shift}
\end{figure}

This will allow us to cluster nodes near a given location $(\bar{r},\bar{\theta})$, for example by using the Bayliss-Turkel map \cite{BaylissTurkel} as $g_1$ and the map introduced by the authors in \cite{BerElef} as $g_2$; see Figure \ref{ChebEquiGRID_Shift} for an example.

Then, considering the interpolant $\calI$ over a grid of $(n_1+1)n_2$ nodes in $[0,2]\times[0,2\pi)$,
 we can reorder the basis functions from $(0,0)$ to $(n_1,n_2-1)$ according to
$$ (0,0)\rightarrow 0, \quad (0,1)\rightarrow 1,\quad  (0,2)\rightarrow 2, \dots, (n_1,n_2-1)\rightarrow n_1(n_2-1) $$
and write the interpolant \eqref{2dBaryInterp} in the linear form
$$\calI f(r,\theta) =\sum_{j=0}^{n_1(n_2-1)} B_{i_j} (r,\theta) f(r_{i_j},\theta_{i_j}) $$
with
$$ B_{i_j} (r,\theta) = b_{i_j}^{(1)} (r) b_{i_j}^{(2)}(\theta) .$$
The Lebesgue constant of $\calI$ is \cite[p.\ 24]{PowellApprox}
$$ \Lambda_{n_1n_2} =\max_{(r,\theta)\in[0,2]\times[0,2\pi)} \sum_{i_j=0}^{n_1(n_2-1)} |B_{i_j}(r,\theta)|, $$
which is equivalent to 
\begin{align*}
    \Lambda_{n_1n_2} &=  \max_{(r,\theta)\in[0,2]\times[0,2\pi]} \sum_{i=0}^{n_1} \sum_{j=0}^{n_2-1} |b_i^{(1)}(r)|\cdot |b_j^{(2)}(\theta)| \\
    &= \bigg( \max_{r\in[0,2]} \sum_{i=0}^{n_1} |b_i^{(1)}(r)| \bigg) \bigg( \max_{\theta\in[0,2\pi]} \sum_{i=0}^{n_2-1} |b_i^{(2)}(\theta)| \bigg) = \Lambda_{n_1}^{(1)} \Lambda_{n_2}^{(2)},
\end{align*}
where $\Lambda_{n_1}^{(1)}$ is the Lebesgue constant of the interpolant $R_1$ for the $n_1+1$ nodes $r_i$ and $\Lambda_{n_2}^{(2)}$ is that of $t_n$ for the $n_2$ nodes $\theta_j$ and with the last equality valid since the variables are separate and the maximum is in a rectangle.

Since it is proven in \cite{LebBerrut2GenNodes} that $\Lambda_{n_1}^{(1)} = \calO(\log (n_1))$ and in \cite{BerElef2} that $\Lambda_{n_2}^{(2)} = \calO(\log (n_2))$, this shows the following.

\begin{thm} \label{Leb2D}
 The Lebesgue constant of the interpolant \eqref{2dBaryInterp} at $n_1+1$ Chebyshev nodes times $n_2$ equispaced nodes as well as that of the interpolant \eqref{2dConfBaryInterp} at conformally shifted nodes obey 
 $$ \Lambda_{n_1 n_2} = \calO (\log n_1 \log n_2).$$
\end{thm}

\section{Interpolation on two-dimensional starlike domains}\label{sec3}

We shall now consider functions on starlike domains \cite[p.\ 92]{Burenkov}. \index{Starlike domains}

\begin{defi}
A domain $\Omega \subset \bbR^2$ is called \emph{starlike with respect to a point} $x\in\Omega$, if for every point $y\in\Omega$ the closed segment $[x,y]$ is contained in $\Omega$. Furthermore, a domain $\Omega\subset \bbR^2$ is called \emph{starlike}, if it is starlike with respect to at least one of its points.
\end{defi}

We shall write $\Omega$ in polar coordinates, with $(0,0)$ the center of the domain, i.e., as the domain contained in the curve
$$ z(\theta) := \Big( \rho(\theta) \cos(\theta),\rho(\theta) \sin(\theta) \Big)$$
for a $2\pi$-periodic function 
$$ \rho(\theta): [0,2\pi] \to \bbR^{> 0},$$
i.e.,
$$ \Omega := \Big\{ (r,\theta)\in\bbR^2 : r<\rho(\theta) \Big\}. $$

\begin{figure}[H]
    \centering
    \graphicspath{ {Images/}}
    \subfigure[]
        {\includegraphics[height=1.74in]{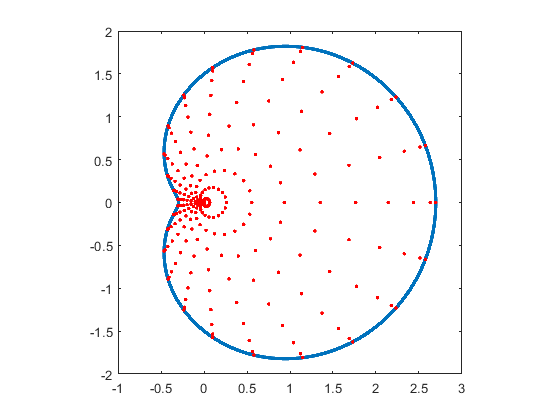}}
     ~ 
    \subfigure[]
        {\includegraphics[height=1.74in]{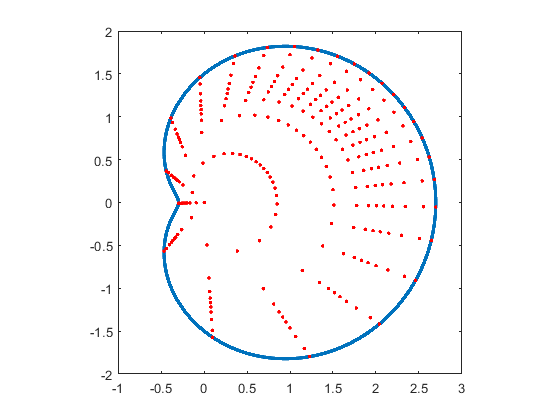}}
    \caption{On the left, an example of the mapped uniform grid, and on the right, an example of a conformally shifted mapped grid.}
    \label{Starlike}
\end{figure}

Our goal is to approximate a function $f$ on $\bar{\Omega}$ by using the interpolant $\calI$ of the former section. 
To that aim, we consider the homothetic grid of points in $\bar{\Omega}$, which consists in the grid (see the example in Figure \ref{Starlike})
$$ G_{i,j} = \Big( r_i \rho(\theta_j), \theta_j \Big) := \Big( \xi_{i,j},\p_j \Big) \in\bOmega,$$
with the corresponding function values 
$$ f_{i,j}=f(\xi_{i,j},\p_j),$$
and we introduce the change of variable from $\bOmega$ to the disk $E$ in \eqref{DiskRad2}
\begin{equation} \label{DomainMap}
    S(\xi,\p)=\left( \frac{2\xi}{\rho(\p)},\p \right)=\Big( S_1(\xi,\p),S_2(\xi,\p) \Big):= ( r,\theta ) .  
\end{equation}
In this way we can interpolate on $E$ using the data of the function on the grid on $\bOmega$ by the function
\begin{equation} \label{DomainBaryInterp}
    \calI [f] (S(\xi,\p)) = \frac{ \sum_{i=0}^{n_1} {\vphantom{\sum}}'' \sum_{j=0}^{n_2-1} \frac{(-1)^{i+j}}{S_1(\xi,\p)-S_1(\xi_{i,j},\p_j)} \cst\left(\frac{S_2(\xi,\p)-S_2(\xi_{i,j},\p_j)}{2} \right) f_{i,j}    }{\sum_{i=0}^{n_1}{\vphantom{\sum}}''\sum_{j=0}^{n_2-1} \frac{(-1)^{i+j}}{S_1(\xi,\p)-S_1(\xi_{i,j},\p_j)} \cst\left(\frac{S_2(\xi,\p)-S_2(\xi_{i,j},\p_j)}{2} \right) }.
\end{equation}

This is similar to a technique introduced by De Marchi et al.\ \cite{FakeNodes1,FakeNodes2,Gibbs} called the Fake Nodes (or mapped basis) Approach, which consists in considering the interpolant constructed on the mapped interpolation nodes and a set of evaluation points, mapped as well.

The interpolant \eqref{DomainBaryInterp} inherits several nice properties from the interpolant $\eqref{2dBaryInterp}$.

\begin{pro} \label{EquivLebesgue}
Let $S$ be as in \eqref{DomainMap}. Then the Lebesgue constant $ \bar{\Lambda}_{n_1 n_2}$ of the interpolant \eqref{DomainBaryInterp} at the set of nodes $G_{i,j}$ satisfies
$$ \bar{\Lambda}_{n_1 n_2} = \calO(\log n_1 \log n_2 ). $$
\end{pro}

\begin{proof}
By definition, the Lebesgue constant of \eqref{DomainBaryInterp} is
\begin{align*}
    \bar{\Lambda}_{n_1 n_2} &= \max_{(\xi,\p)\in\bOmega} \sum_{i_j=0}^{n_1(n_2-1)} |B_{i_j}(S(\xi,\p))| = \max_{(r,\theta)\in S(\bOmega)} \sum_{i_j=0}^{n_1(n_2-1)} |B_{i_j}(r,\theta)| \\
    &=\max_{(r,\theta)\in E} \sum_{i_j=0}^{n_1(n_2-1)} |B_{i_j}(r,\theta)| = \Lambda_{n_1 n_2},
\end{align*} 
where $\Lambda_{n_1 n_2}$ is the Lebesgue constant of the interpolant \eqref{2dBaryInterp}.
Theorem \ref{Leb2D} then yields the result.
\end{proof}

Note that this is true in particular when the points $G_{i,j}$ are such that their images under the map $S$ are conformally shifted nodes in the disk.

Moreover, we may consider the effect of perturbed data $\widetilde{f}(G_{i,j})$ on the interpolant.

\begin{pro}
Let $S$ and the interpolant $\calI[f]$ be as above and consider perturbed data $\widetilde{f}$. Then,
$$ \max_{(\xi,\p)\in\Omega} |\calI[f](\xi,\p) - \calI[\widetilde{f}](\xi,\p)| \leq \Lambda_{n_1 n_2}  \, \max_{(\xi,\p)\in\Omega} |f(\xi,\p)-\widetilde{f}(\xi,\p)|. $$
\end{pro}

\begin{proof}
This follows directly from the classical estimate on the conditioning of an interpolant (see e.g. \cite[p.\ 23]{PowellApprox}), i.e., 
$$ \norm{p_n[f]-p_n[\widetilde{f}]} \leq \Lambda_n \norm{f-\widetilde{f}}, $$
with $p_n[h]$ the polynomial of degree at most $n$ which interpolates the function $h$ at a fixed set of nodes and $\Lambda_n$ the norm of the interpolation operator.
\end{proof}

\section{Convergence}

In order to study its convergence in a starlike domain, we first consider the interpolant in the disk (or the rectangle).
Then, we can prove the following bound.

\begin{thm} \label{TensorError} 
Let $f\in\calC^\infty(B)$ be such that for any fixed $r$ the function $f(r,\cdot)$ satisfies the hypothesis of Theorem \ref{TrigoError} and for any fixed $\theta$ the function $f(\cdot,\theta)$ satisfies the hypothesis of Theorem \ref{ChebError}. Then the error of the interpolant \eqref{2dBaryInterp} can be bounded as
\begin{equation} \label{errbound1}
\norm{ f - \calI [f] } \leq \Big( 2\Lambda_{n_1}^{(1)} + \Lambda_{n_1}^{(1)}\Lambda_{n_2}^{(2)} +1 \Big) \max\left\{ \frac{4M_1}{\sigma-1}\sigma^{-n_1}, 2 M_2 \ctg\left( \frac{a}{2}\right) e^{-a n_2} \right\} .
\end{equation}
\end{thm}

\begin{proof}
First, let us consider the points $(\bar{r},\bar{\theta})\in B$ such that
$$ (\bar{r},\bar{\theta})=\arg\max_{(r,\theta)\in B} |f(r,\theta)-\calI[f](r,\theta)|. $$
Then, with the notations 
$$ \calI_1[f](r,\theta):=\sum_{i=0}^{n_1} b_i^{(1)}(r) f(r_i,\theta) \qquad \text{and}\qquad \calI_2[f](r,\theta):=\sum_{j=0}^{n_2-1} b_j^{(2)}(\theta) f(r,\theta_j),$$
we can compute
\begin{align*}
    \norm{f-\calI[f]} &= \Big| f(\bar{r},\bar{\theta})-\calI[f](\bar{r},\bar{\theta}) \Big| = \Big| f(\bar{r},\bar{\theta}) - \calI_1 [\calI_2 [f]](\bar{r},\bar{\theta}) \Big| \\
    &= \Big| f(\bar{r},\bar{\theta}) -\calI_1[f](\bar{r},\bar{\theta}) +\calI_1[f](\bar{r},\bar{\theta}) - \calI_1 [\calI_2 [f]](\bar{r},\bar{\theta}) \Big| \\
    &\leq \Big| f(\bar{r},\bar{\theta}) -\calI_1[f](\bar{r},\bar{\theta}) \Big|  +\Big| \calI_1[f](\bar{r},\bar{\theta}) - \calI_1 [\calI_2 [f]](\bar{r},\bar{\theta}) \Big| \\
    &= \Big| f(\bar{r},\bar{\theta}) -\calI_1[f](\bar{r},\bar{\theta}) \Big|  +\Big| \calI_1\big[ f(\cdot,\bar{\theta}) -  \calI_2 [f](\cdot,\bar{\theta})\big](\bar{r}) \Big|,
\end{align*}
where the last line holds in view of the the linearity of $\calI_1$. Thus,
\begin{align*}
    \norm{f-\calI[f]} &\leq \Big| f(\bar{r},\bar{\theta}) -\calI_1[f](\bar{r},\bar{\theta}) \Big|  + \norm{\calI_1} \Big| f(\bar{r},\bar{\theta}) -  \calI_2 [f](\bar{r},\bar{\theta}) \Big| \\
    &\leq (1+\Lambda_{n_1}^{(1)}) \frac{4M_1}{\sigma-1} \sigma^{-n_1} + \Lambda_{n_1}^{(1)}(1+\Lambda_{n_2}^{(2)}) 2M_2 \ctg\left(\frac{a}{2}\right) e^{-a n_2}, 
\end{align*}
where we have used the classic estimate for interpolation \cite[p.\ 24]{PowellApprox}
$$ \norm{h-p_n[h]}\leq (1+ \Lambda_n) \norm{h-p^*}, $$
where $p^*$ denotes the best approximation of $h$ in $\bbP_n$ and the fact that, since $\calI_1$ (respectively $\calI_2$) corresponds to the interpolation at Chebyshev nodes of the second kind (resp. trigonometric interpolation at equidistant nodes) we have from Theorem \ref{ChebError} that, for a function $h:\bbR\to\bbR$ analytic in the ellipse $\calE_\sigma$,
$$ \norm{h-p^*} \leq \norm{h-\calI_1[h]}\leq \frac{4M_1}{\sigma-1}\sigma^{-n_1},$$
and similarly for $\calI_2$ with Theorem \ref{TrigoError}.

Therefore, we have \eqref{errbound1}.
\end{proof}

Similarly we can prove a convergence theorem with weaker conditions, by proceeding as in the proof of Theorem \ref{TensorError}.

\begin{thm}
Let $f$ be a function such that for any fixed $r$ the function $f(r,\cdot)$ satisfies the hypothesis of Theorem \ref{TrigoErrorJump} and for any fixed $\theta$ the function $f(\cdot,\theta)$ satisfies the hypothesis of Theorem \ref{ChebyErrorJump}. Then the error of the interpolant \eqref{2dBaryInterp} can be bounded as
\begin{equation}
    \norm{f-\calI[f]}\leq \Big( 2\Lambda_{n_1}^{(1)} + \Lambda_{n_1}^{(1)}\Lambda_{n_2}^{(2)} +1 \Big) \max\left\{ \frac{4V_1}{\pi \nu (n-\nu)^{\nu}}, \frac{V_2}{\pi n^\mu}\left(\frac{1}{n}+\frac{2}{\mu} \right) \right\} .
\end{equation}
\end{thm}

Finally, we obtain similar bounds for the interpolant at conformally shifted nodes, by following the proof of Theorem \ref{TensorError} but estimating with the best approximation in the spaces  $\calR_n$ and $\bbSigma_n$ introduced in section \ref{Sec2}, instead of the best approximation in the polynomial spaces.

\begin{thm} \label{TensorConfError}
Let $f\in\calC^\infty(B)$ be such that for any fixed $r$ the function $f(r,\cdot)$ satisfies the hypothesis of Theorem \ref{BalThmConf} and for any fixed $\theta$ the function $f(\cdot,\theta)$ satisfies the hypothesis of Theorem \ref{TransfConv}. Then, we can bound the error of the interpolant \eqref{2dConfBaryInterp} as
\begin{equation} 
\norm{ f - \calI [f] } \leq \Big( 2\Lambda_{n_1}^{(1)} + \Lambda_{n_1}^{(1)}\Lambda_{n_2}^{(2)} +1 \Big) \max\left\{ C_1 \sigma^{-n_1}, C_2 e^{-a n_2} \right\},
\end{equation}
for two constants $C_1,C_2$, where $C_1$ depends on $\sigma$ and $g_1$, $C_2$ on $a$ and $g_2$, and both on the function $f$.
\end{thm}

Finally, note that the interpolant in the starlike domain inherits the error from the interpolant on the disk.

\begin{thm}
Let $S$ be as in \eqref{DomainMap}, $\calI$ as in \eqref{2dBaryInterp} (or \eqref{2dConfBaryInterp}), and $\calI^S$ as in $\eqref{DomainBaryInterp}$. Then for a function $f:\Omega\to\bbR$ we have
\begin{equation}
    \norm{f-\calI^S[f]}_{\Omega} = \norm{h-\calI[h]}_{E},
\end{equation}
where $h=f\circ S^{-1}$.
\end{thm}

\begin{proof}
By definition we have $\calI^S=\calI \circ S$; then, if we consider as $h$ the function such that $h\circ S=f$, we have
$$ \norm{\calI^S[f]-f}_{\Omega} =\norm{\calI\circ S- h\circ S}_{\Omega} = \norm{\calI[h]-h}_{S(\Omega)}, $$
and $S(\Omega)=E$.
\end{proof}

Consequently, since the interpolant converges as that of the function $h=f\circ S^{-1}$, where 
$$ S^{-1}(r,\theta) := \left( \frac{r\rho(\theta)}{2},\theta \right), $$
the convergence depends on $f$ and $S$; a smoother function $\rho$ will lead to a better interpolant \eqref{DomainBaryInterp} of the function $f$ in $\bOmega$.

We have also interpolated in non-smooth domains but the results were not satisfactory, as the function $h$ itself is non-smooth along the boundary, i.e., for $r=2$.
A way of dealing with such a non-smooth function $\rho$, as for example the polar representation of the square
\begin{equation} \label{RhoSquare}
    \rho(\theta) = \min \left\{ \frac{1}{|\cos(\theta)|}, \frac{1}{|\sin(\theta)|}  \right\},
\end{equation} 
is to replace it with a close enough smoother approximant $\widetilde{\rho}$ and interpolate in the approximate domain $\widetilde{\Omega} = \widetilde{S}(E)$ instead of $\bOmega=S(E)$, where $\widetilde{S}$ is the function $S$ with $\widetilde{\rho}$ in place of $\rho$. 
Some tests in this direction are reported in the next section.

\section{Numerical tests}

In this section we test the interpolant \eqref{DomainBaryInterp} on domains inside the curve 
$$ z(\theta) := \Big( \rho(\theta) \cos(\theta),\rho(\theta) \sin(\theta) \Big),$$
for a $2\pi$-periodic boundary parametrization
$$ \rho(\theta): [0,2\pi] \to \bbR^{> 0}.$$

\begin{figure}
\centering
\graphicspath{ {Images/}}
    \subfigure[Lima\c{c}on]
        {\includegraphics[height=1.74in]{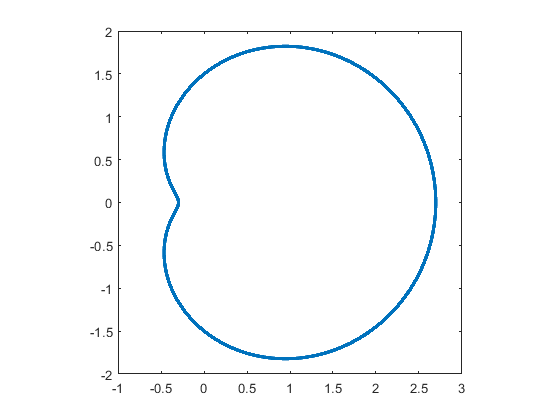}}
     ~ 
    \subfigure[First butterfly]
        {\includegraphics[height=1.74in]{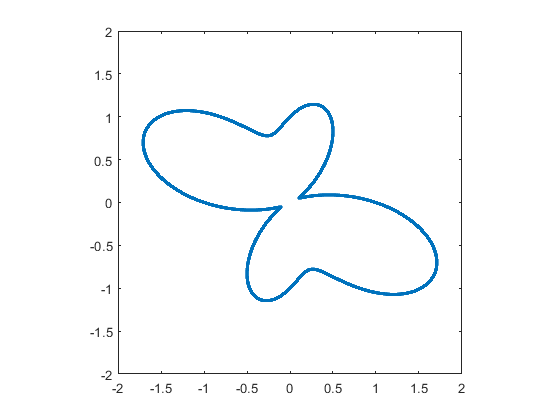}}
     ~ 
    \subfigure[Second butterfly]
        {\includegraphics[height=1.74in]{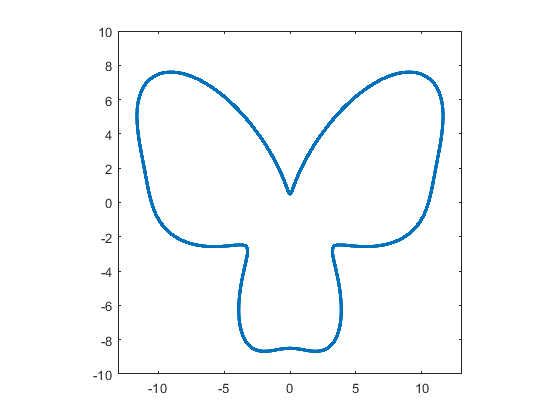}}
     ~ 
    \subfigure[Asterisk]
        {\includegraphics[height=1.74in]{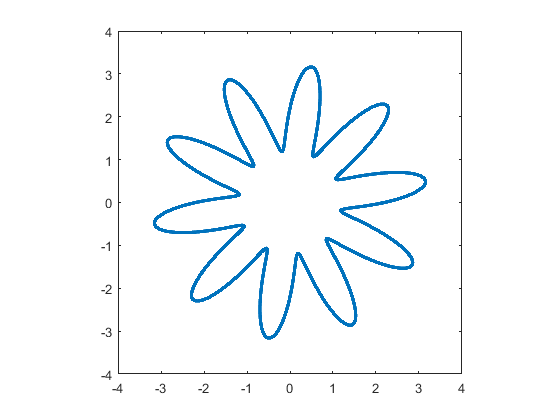}}
     ~ 
         \caption{Test domains} 
    \label{Domains}
\end{figure}

We first consider the following domains (see Figure \ref{Domains}):

\begin{itemize}
    \item The lima\c{c}on
     $$ \rho_1(\theta) = 1.5+1.2\cos(\theta);$$
     \item A first butterfly-shaped domain
     $$ \rho_2(\theta) = 1-\cos(\theta)\sin(3\theta); $$
     \item A second butterfly-shaped domain
     $$ \rho_3(\theta) = 7.5-\sin(\theta)+4\sin(3\theta)-\sin(7\theta)+3\cos(2\theta);$$
     \item The asterisk
     $$ \rho_4(\theta) = \sin(10\,\theta)+2.2.$$
\end{itemize}

To estimate the error, we consider a grid of $170\times 170$ uniformly spaced points in a rectangle around the domain (see Table \ref{Box} for the rectangles) and we filter them by considering only the points inside the domain. Then, we compute the maximum of the absolute value of the difference between the interpolant and the function in these points.

\begin{table}
\centering
\begin{tabular}{l||c}
 & Rectangle  \\
\hline
$\rho_1$  & $[-1, 3]\times[-2, 2]$   \\
$\rho_2$  & $[-2, 2]\times [-2, 2]$  \\
$\rho_3$  & $[-13, 13]\times[-10, 10]$   \\
$\rho_4$  & $[-4, 4]\times [-4, 4]$   \\
$\rho_5$  & $[-2, 2]\times [-2, 2]$  
\end{tabular}
\caption{Rectangles for the evaluation points with the various domains}
\label{Box}
\end{table}

Notice that all these domains are constructed via a smooth function $\rho$; thus the expected convergence behaviour will not be influenced by infinitely smooth changes of variable such as conformal point shifts.

\begin{figure}
    \centering
    \graphicspath{ {Images/}}
    \subfigure[]
        {\includegraphics[height=1.74in]{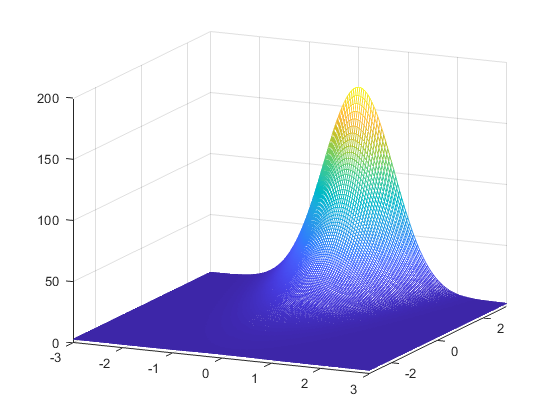}}
     ~ 
    \subfigure[]
        {\includegraphics[height=1.74in]{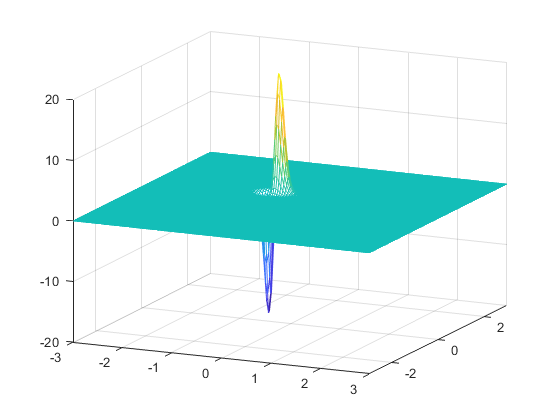}}
    \caption{On the left, the first test function $f_1$, and on the right, the second test function $f_2$ with the parameters used in the computations.}
    \label{TestFun}
\end{figure}

We document our computations with two infinitely smooth functions, a first one which does not require point shifts,
$$ f_1(x,y) = 3e^{-x^2+y+1}+3 ,$$
(on the left in Figure \ref{TestFun}) and for which we expect (approximate) exponential convergence.

As we observe in Table \ref{Er_f1}, each error is approximately the square of that on the previous row when $(n_1,n_2)$ is doubled; this indeed reflects an exponential decay of the error.

\begin{table}
\centering
\begin{tabular}{l||c|c|c|c}
$(n_1,n_2)$ & $\rho_1$ & $\rho_2$ & $\rho_3$ & $\rho_4$ \\
\hline
(10,30)   & 1.6762e--02 & 1.3439e--01 & 1.4178e+01 & 2.8832e+01  \\
(20,60)   & 1.6080e--07 & 3.3468e--04 & 2.1093e+00 & 3.0920e+00  \\
(40,120)  & 8.5265e--14 & 1.3499e--10 & 9.0279e--02 & 1.5704e--02  \\
(80,240)  & 1.2790e--13 & 7.1054e--14 & 2.0515e--05 & 4.6051e--07  \\
(160,480) & 1.4921e--13 & 1.0303e--13 & 9.9476e--14 & 5.6843e--13 
\end{tabular}
\caption{Errors with the function $f_1$ in the various domains}
\label{Er_f1}
\end{table}

To take advantage of the ability of linear barycentric rational interpolation to accommodate steep gradients (fronts) by simply replacing the nodes in the "classical" interpolant with shifted ones which accumulate in the vicinity of the fronts, i.e., $f_{ij}\coloneqq f(g_1(r_{i}) \rho(g_2(\theta_j)),g_2(\theta_j))$, we use the interpolant $\calI[f]$ in \eqref{DomainBaryInterp} and we consider secondly a function (on the right in Figure \ref{TestFun}) with a front in a precise location, i.e.,
$$ f_2(x,y) = 40 \frac{\erf\left( \sqrt{\epsilon/2}(x+0.6) \right)}{\erf(\sqrt{\epsilon/2})} e^{-30(x+0.6)^2} e^{-60(y-0.6)^2},$$
where $\epsilon=100$ and $\erf$ is the error function; $f_2$ has a front in $(-0.6,0.6)$ (which becomes $(0.6\sqrt{2},\frac{7\pi}{4})$ in polar coordinates). We use the conformal maps 
$$g_1 =\beta + \frac{1}{\alpha} \tan\left(  \lambda (x-\mu) \right), \quad \text{ and }\quad g_2(\theta) = - i \log \left( \frac{e^{i \p}+\eta e^{i\theta}}{1+e^{i\p} {\eta e^{-i\theta}}} \right)$$
with $\beta$ and $\p$ such that the nodes cluster around the location of the front and  for the density $\alpha = 2.8$ for $g_1$  and $\eta = 0.65$ for $g_2$; $g_1$ and $g_2$ are respectively the Bayliss-Turkel map \cite{BaylissTurkel} and the map introduced by the authors in \cite{BerElef} (since everything is smooth, the error is horizontal at the optimal values and the results are not very sensitive to the parameters there, see the second column of the tables in \cite{BerMitPROC}).
The results, displayed in Tables \ref{Er_f2_1} and \ref{Er_f2_2}, indeed show that, as in the unidimensional case, we achieve a much faster convergence by using conformal maps to cluster nodes in the vicinity of the location of a front instead of the plain Chebyshev and equispaced points.

\begin{table}
\centering
\begin{tabular}{l||c|c|c|c}
$(n_1,n_2)$ & $\rho_1$ & $\rho_1$ conf & $\rho_2$ & $\rho_2$ conf \\
\hline
(10,30)   & 2.2524e+01 & 4.9054e+00 & 1.7898e+01 & 1.8408e+00  \\
(20,60)   & 7.7530e+00 & 1.7487e--02 & 4.6606e+00 & 3.7443e--02  \\
(40,120)  & 1.0473e--01 & 6.2046e--07 & 6.1903e--02 & 1.0631e--05  \\
(80,240)  & 2.2811e--07 & 1.8474e--13 & 1.5352e--06 & 5.8037e--13 
\end{tabular}
\caption{Errors with the function $f_2$ for the first two domains}
\label{Er_f2_1}
\end{table}

\begin{table}
\centering
\begin{tabular}{l||c|c|c|c}
$(n_1,n_2)$ & $\rho_3$ & $\rho_3$ conf & $\rho_4$ & $\rho_4$ conf \\
\hline
(10,30)   & 1.8077e+01 & 1.3739e+01 & 2.2392e+01 & 1.3262e+01  \\
(20,60)   & 7.6290e+00 & 2.7313e+00 & 2.1117e+01 & 5.3838e+00  \\
(40,120)  & 1.3786e+00 & 2.6799e--02 & 1.1580e+01 & 7.5581e--01  \\
(80,240)  & 1.9880e--02 & 1.3075e--06 & 9.3368e--01 & 7.3685e--03  \\
(160,480) & 2.3293e--08 & 1.0303e--13 & 1.5659e--03 & 1.3545e--08 
\end{tabular}
\caption{Errors with the function $f_2$ for the next two domains}
\label{Er_f2_2}
\end{table}

\vspace{6pt}
Furthermore, we consider a domain which corresponds to a non-smooth $\rho$, the square $[-1,1]^2$, parametrized by the function \eqref{RhoSquare} which we denote by $\rho_5$.
In this case, approximating functions with the transplanted interpolant often leads to catastrophic results, see the columns $\rho_5$ in Table \ref{FakeRho}. To avoid this, we slightly modify the domain by approximating $\rho_5$ with a smooth interpolant $\widetilde{\rho_5}$ and we use the interpolant transplanted via the function
$$ \tilde{S}(\xi,\p)=\left( \frac{2\xi}{\widetilde{\rho_5}(\p)},\p \right)=\Big( \widetilde{S}_1(\xi,\p),\widetilde{S}_2(\xi,\p) \Big):= ( r,\theta ).  $$

For Table \ref{FakeRho} we approximate $\rho_5$ via the approximant produced by the AAA algorithm \cite{AAA}: the corresponding $\widetilde{\rho}_5$ clearly gives a better interpolant than $\rho_5$. Unfortunately, because of its "almost corner", the AAA approximant may lead, in practice, to an interpolant with a slower than exponential convergence in the neighbourhood of the corners; furthermore, the interpolant is obtained merely in an approximation of the original domain $[-1,1]^2$.   

\begin{table} 
\centering
\begin{tabular}{l||c|c||c|c|c}
 & \multicolumn{1}{r}{$f_1$}  & & \multicolumn{3}{c}{$f_2$} \\
\hline
$(n_1,n_2)$ & $\rho_5$ & $\widetilde{\rho}_5$ & $\rho_5$ & $\widetilde{\rho}_5$ & $\widetilde{\rho}_5$ conf \\
\hline
(10,30)   & 2.9917e+00   & 4.4789e--01 & 9.3838e+03   & 1.2996e+01 & 1.8835e+00  \\
(20,60)   & 3.1510e+00   & 2.6646e--01 & 1.5954e+09   & 7.4563e+00 & 1.3315e+00  \\
(40,120)  & \textrm{Inf} & 7.7975e--02 & \textrm{Inf} & 1.9865e+00 & 1.0297e--01  \\
(80,240)  & \textrm{Inf} & 3.7900e--02 & \textrm{Inf} & 4.9346e--01 & 4.0428e--02  \\
(160,480) & \textrm{Inf} & 4.8907e--03 & \textrm{Inf} & 2.1970e--01 & 8.1943e--03 
\end{tabular}
\label{FakeRho}
\caption{Errors with the $\rho_5$ and its approximation $\widetilde{\rho_5}$}
\end{table}

\begin{figure}
\centering
\graphicspath{ {Images/}}
    \subfigure[]
        {\includegraphics[height=1.74in]{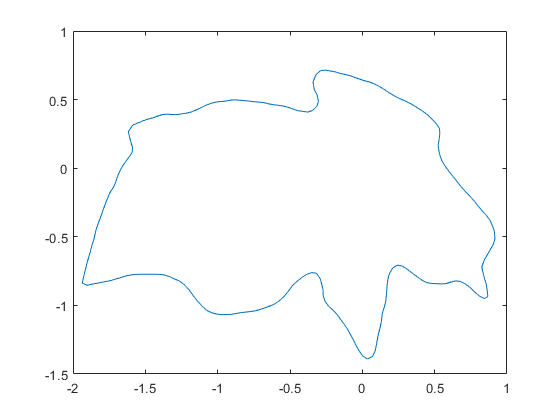}}
     ~ 
    \subfigure[]
        {\includegraphics[height=1.74in]{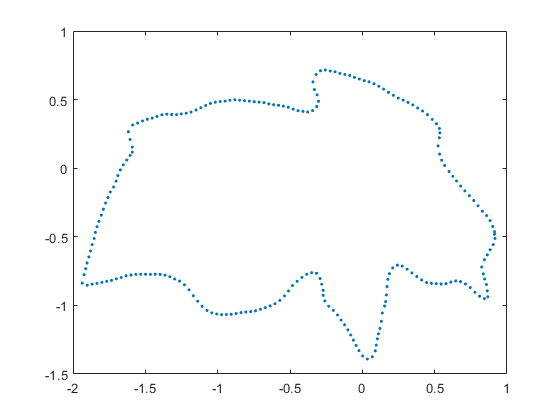}}
     ~ 
         \caption{Swiss-like domain: (a) The sketch; (b) The extracted nodes} 
    \label{SwisslikeNodes}
\end{figure}

\begin{figure}
\centering
\graphicspath{ {Images/}}
    \subfigure[]
        {\includegraphics[height=1.74in]{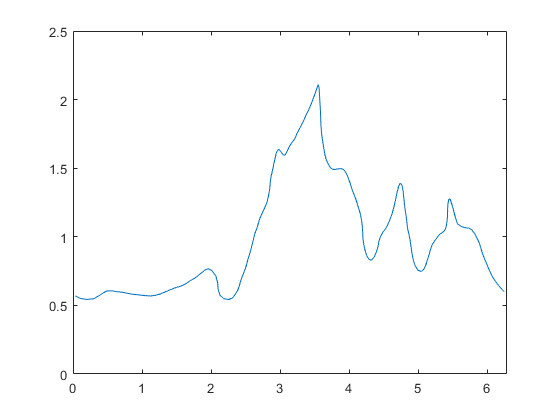}}
     ~ 
    \subfigure[]
        {\includegraphics[height=1.74in]{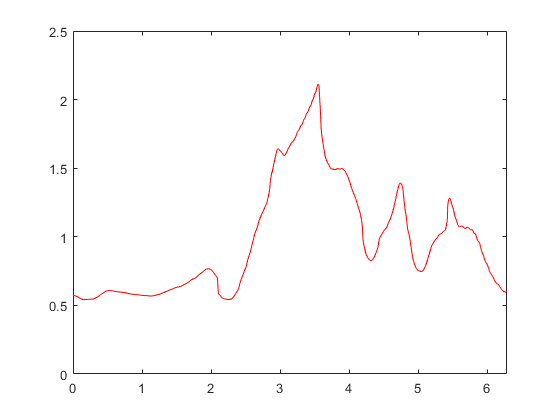}}
     ~ 
         \caption{(a) $\rho_6$; (b) The approximation $\widetilde{\rho}_6$} 
    \label{SwisslikeRHO}
\end{figure}

\vspace{6pt}
As a last test, we consider a sketch of Switzerland, in which we simplify the boundary in order to have a domain which is starlike with respect to the origin (see Figure \ref{SwisslikeNodes}). To obtain the boundary curve, we extracted some points which we use as nodes to interpolate the unknown $\rho_6$ representing the boundary of the "Swiss-like" domain; then we approximate the curve with Berrut's first interpolant $R_0$ (see e.g. \cite{BerInt1,BaryHormann}), which produces an approximation, let us denote it by $\widetilde{\rho}_6$, of the unknown $\rho_6$ (see Figure \ref{SwisslikeRHO}). In this case, the rectangle in which we compute the error of the interpolant is $[-2.5, 1.5]\times [-2, 1.5]$. In Table \ref{FakeRhoSwiss} we display the errors produced with the functions $f_1$ and $f_2$; the exponential convergence does not yet show up without conformal shift, but it does with such a shift.

\begin{table} 
\centering
\begin{tabular}{l||c||c|c}
 & $f_1$  & \multicolumn{2}{c}{$f_2$} \\
\hline
$(n_1,n_2)$ &  $\widetilde{\rho}_6$ &  $\widetilde{\rho}_6$ & $\widetilde{\rho}_6$ conf \\
\hline
(40,120)    & 4.3451e--01 & 9.6853e+00 & 5.3945e+00  \\
(80,240)    & 3.2449e--01 & 9.2768e+00 & 6.3695e--01  \\
(160,480)   & 1.5968e--01 & 6.2080e+00 & 3.0078e--01  \\
(320,960)   & 4.1611e--02 & 1.3068e+00 & 2.6339e--02 \\
(640,1920)  & 1.9271e--02 & 4.1124e--01 & 8.5433e--05 \\
(1280,3840) & 1.7645e--03 & 5.4788e--02 & 1.0787e--09
\end{tabular}
\caption{Errors for the "Swiss-like" domain}
\label{FakeRhoSwiss}
\end{table}

\section{Concluding remarks}

The present work has introduced a simple generalisation of linear barycentric
rational interpolants to two-dimensional domains. The idea is to parameterize a
(starlike) domain with polar coordinates and to use linear barycentric interpolation
in the radial direction and its trigonometric version in the circular
(homothetic) direction. (In reality, the interpolation happens in a disk, to which the
original problem is transplanted.)
Up to a logarithmic factor (which arises from the proof, but may not appear in practice),
the resulting tensor-product-like interpolant converges exponentially when the function
is analytic, and as $\calO(h^\nu)$ when $f^{(\nu)}$ has bounded variation.
Impressive numerical examples amply confirm these theoretical convergence results.

We do not want to conceal certain limitations of the method. The most important
seems to be the fact that the rapid convergence requires that the boundary
parametrization
is as smooth as the interpolated function $f$, as the interpolant involves the
trigonometric interpolation of $f$ along the boundary.
Another one
is the fact that the boundary curve cannot be arbitrary: for instance, its interior must
contain a point, to be chosen as the center of the domain, from which the
representation of the curve does not have too large a derivative.

One can easily generalize the method to accomodate functions with several fronts by the method introduced in \cite{BerMit1} and extended to the circular case in \cite{BerElef}. But the required computing power could rapidly approach the limitations of MATLAB on a personal computer.

Another issue is that the point shifts deplete the grid close to the boundary (if the steep fronts are not there). In case $f$ has a disturbing steep gradient outside $\Omega$, one can ensure that a radius passes through
the corresponding abscissa and attach poles to the radial interpolant, symmetrically with respect to the radius; that way one may expect to gain about 2–3 digits \cite{BerMit1}; the barycentric
representation is ideal for this attachment.

Finally, we note that
one nice feature of this interpolant is the simplicity of the formulae for its partial
derivatives along the lines making up its grid: in radial direction they are given
by Schneider and Werner's formula \cite{SchneiderWerner},
for the circular direction one finds them in \cite{Baltensperger} (up to
inner derivatives of the variable transformation in the chain rule).
These formulae could be used to approach the (infinitely smooth) solutions of partial differential equations such as Poisson problems in starlike domains. However, the fact that the mesh is not orthogonal will lead to more complicated formulae for the Laplacian (chain rule) than just those of Schneider-Werner and Baltensperger.

\section*{Acknowledgements}
The authors thank the referees for their careful reading
of the manuscript
and their numerous suggestions, which have improved this work.
The research of the second author has been performed within the Rete ITaliana di Approssimazione (RITA), the UMI Group TAA “Approximation Theory and Applications”, and with the support of GNCS-IN$\delta$AM.
\vspace{14pt}

\printbibliography

\end{document}